\documentclass[bsl]{asl}
\title{Grothendieck rings of $\Z$-valued fields }
\author{Raf Cluckers$^*$}
\revauthor{Cluckers, Raf}
\address{Department of Mathematics\\
Katholieke Universiteit Leuven\\ Celestijnenlaan 200B\\ B-3001 Heverlee,
Belgium}
 \email{raf.cluckers@wis.kuleuven.ac.be}
\author{Deirdre Haskell}
\revauthor{Haskell, Deirdre}

\thanks{$^*$ Research Assistant of the Fund for Scientific Research --
 Flanders (Belgium)(F.W.O.)}

\address{Department of Mathematics and Statistics\\
McMaster University\\
1280 Main St. West\\
Hamilton, Ontario, Canada\\
L8S 4K1}

\email{haskell@math.mcmaster.ca}

\newtheorem{theorem}{Theorem}
\newtheorem{proposition}{Proposition}

\newtheorem{cor}{Corollary}
\theoremstyle{definition}
\newtheorem{definition}{Definition}
\newtheorem*{remark}{Remark}

\newcommand{\Q}{\ensuremath{\mathbb{Q}}}
\newcommand{\Z}{\ensuremath{\mathbb{Z}}}
\newcommand{\N}{\ensuremath{\mathbb{N}}}

\newcommand{\Lm}{\ensuremath{\mathcal{L}}}
\newcommand{\M}{\ensuremath{\mathcal{M}}}
\newcommand{\F}{\ensuremath{\mathbb{F}}}
\newcommand{\Def}{\ensuremath{\mathcal{D}ef}}
\begin{document}

\begin{abstract}
We prove the triviality of the Grothendieck ring of a $\Z$-valued field
$K$ under slight conditions on the logical language and on $K$. We
construct a definable bijection from the plane $K^2$ to itself minus a
point. When we specialize to local fields with finite residue field, we
construct a definable bijection from the valuation ring to itself minus a
point.
\end{abstract}

\maketitle

At the Edinburgh meeting on the model theory of valued fields in
May 1999, Luc B\'elair posed the question of whether there is a
definable bijection between the set of $p$-adic integers and the
set of $p$-adic integers with one point removed. At the same
meeting, Jan Denef asked what is the Grothendieck ring of the
$p$-adic numbers, as did Jan Kraj\'{\i}\v{c}ek independently in
\cite{K}. A general introduction to Grothendieck rings of logical
structures was recently given in \cite{KS} and in [DL2, par. 3.7].
Calculations of non-trivial Grothendieck rings and related topics
such as motivic integration can be found in \cite{DL} and
\cite{DL2}. The logical notion of the Grothendieck ring of a
structure is analogous to that of the Grothendieck ring in the
context of algebraic $K$-theory and has analogous elementary
properties (see \cite{S}). Here we recall the definition.

\begin{definition} Let $\M$ be a structure and $\Def(\M)$ the set
of definable subsets of $M^n$ for every positive integer $n$. For any
$X,Y\in\Def(\M)$, write $X\cong Y$ iff there is a definable bijection (an
isomorphism) from $X$ to $Y$. Let $F$ be the free abelian group whose
generators are isomorphism classes $\lfloor X\rfloor$ with $X\in\Def(\M)$
(so $\lfloor X\rfloor =\lfloor Y\rfloor$ if and only if $X\cong Y$) and
let $E$ be the subgroup generated by all expressions $\lfloor X\rfloor
+\lfloor Y\rfloor -\lfloor X\cup Y\rfloor -\lfloor X\cap Y\rfloor$ with
$X,Y\in\Def(\M)$. Then the Grothendieck group of $\M$ is the quotient
group $F/E$. Write $[X]$ for the image of $X\in\Def(\M)$ in $F/E$. The
Grothendieck group has a natural structure as a ring with multiplication
induced by $[X]\cdot[Y]=[X\times Y]$ for $X,Y\in\Def(\M)$. We call this
ring the Grothendieck ring $K_0(\M)$ of $\M$.
\end{definition}

\noindent It is easy to see that the above questions are related: the
Grothendieck ring is trivial if and only if there is a definable bijection
between $M^k$ and itself minus a point for some $k$, which happens if and
only if the
Grothendieck group is trivial. Moreover, if we find such a $k$ then we
have for any $X\in\Def(\M)$ a definable bijection from the disjoint union
of $M^k\times X$ and $X$ to $M^k\times X$; if there is a definable
injection from $M^k$ into $X$ we find a definable bijection from $X$ to
itself minus a point.

In this paper we answer the questions posed by B\'elair and Denef.
Furthermore, we prove the triviality of the Grothendieck ring of any
$\Z$-valued field which satisfies some slight conditions and give in this
general setting an explicit bijection from the plane to itself minus a
point. For the fields $\Q_p$ and $\F_q((t))$ we explicitly construct a
definable bijection from the valuation ring to itself minus a point.

Dave Marker independently produced a definable bijection from $\Z_p$ to
$\Z_p\setminus\{0\}$, after it was noticed by Lou~van den Dries that its
existence followed from unpublished notes of the second author. The first
author has proved further that there is a definable bijection between any
two definable sets in the $p$-adics if and only if they have the same
dimension. This will appear in a later paper. We thank the referee for
encouraging us to present these results in greater generality than had
been our original intention.

Fix a $\Z$-valued field $K$, that is,  a field with a valuation
$v:K^\times\to Z$ to an ordered group $Z$ which is elementarily
equivalent to the integers in the Presburger language. Let
$R=\{x\in K|v(x)\geq0\}$ be the valuation ring,
$R^*=R\setminus\{0\}$ and $\bar K=R/m$ the residue field, with $m$
the maximal ideal of $R$ and natural projection $R\to \bar
K:x\to\bar x$. An angular component map is a homomorphism
$ac:K^\times\to\bar K^\times$ such that $ac(x)=\bar x$ if
$v(x)=0$. We extend $ac$ to a map $ac:K\to\bar K$ by putting
$ac(0)=0$ (for the existence of angular component maps, see
\cite{P} and \cite{B}).

\begin{definition} Let $\Lm$ be an extension of the
language of rings with $K$ as a model. We say that the structure
$(K,\Lm)$ satisfies condition ($*$) if we can choose an angular component
map $ac$ and an $\Lm$-definable element $\pi\in R$ with $v(\pi)=1$ and
$ac(\pi)=1$ such that the sets $R$ and $R^{(1)}=\{x\in R|ac(x)=1\}$ are
$\Lm$-definable.
\end{definition}

Notice that if condition ($*$) is satisfied, the set $\{(x,y)\in
K^2|v(x)\leq v(y)\}$ is $\Lm$-definable by the formula $\exists z\in R\,(zx=y)$.
A bijection $X\to Y$ with $X,Y\in\Def{(K,\Lm)}$ with $\Lm$-definable graph
will be called an isomorphism.

Let $X\subset K^m$ and $Y\subset K^n$ be definable sets, $m\geq n$. Let
$X'=\{0\}\times X$ and $Y'=\{1\}^{m-n+1}\times Y$. Then we define the
disjoint union $X\sqcup Y$ of $X$ and $Y$ up to isomorphism to be $X'\cup
Y'$.  We say that a set $W$ is isomorphic to $X\sqcup Y$ if $W$ is
isomorphic to $X'\cup Y'$ and then obviously $[W]=[X]+[Y]$.  If $(K,\Lm)$
satisfies condition ($*$) then we can find $W\subset R^m$ with $W\cong
X\sqcup Y$ as follows. The map $i:K\to R$ which sends $x$ to $\pi x$ if
$v(x)\geq0$ and to $1+1/x$ if $v(x)<0$ is a definable injection. For
$m=n=1$, put $X''=\pi.i(X)$ and $Y''=1+\pi.i(Y)$. Then $X''\cong X$,
$Y''\cong Y$ and $X''\cap Y''=\phi$, so  $W=X''\cup Y''$ is isomorphic to
$X\sqcup Y$. For $m>1$, use the same method in each coordinate.

\begin{proposition}\label{bijections} Let $K$ be a $\Z$-valued field, which
is a model for the language $\Lm$. If the structure $(K,\Lm)$ satisfies
condition ($*$), then the following holds:
\item{(i)} The disjoint union of $R$ and $R^{(1)}$ is isomorphic to
$R^{(1)}$ and thus $[R]=0.$
\item{(ii)} The disjoint union of two copies of $R^{* 2}$ is isomorphic to
$R^{* 2}$ itself, and hence $[R^{* 2}]=0.$
\end{proposition}
\begin{proof}
(i) The map
\[
 \{0\}\times R\cup \{1\}\times R^{(1)}\to R^{(1)}:\left\{\begin{array}{rcl}
              (0,x) & \mapsto & 1+\pi x,\\
              (1,x) & \mapsto & \pi x,\end{array}\right.\
\]
is easily seen to be an isomorphism as required. This yields in the
Grothendieck ring $[R] + [R^{(1)}]= [R^{(1)}]$, so $[R]=0$.

(ii)
 Define the sets
\begin{eqnarray*}
X_1=\{(x,y)\in R^{* 2}|v(x)\leq v(y)\},\\
X_2=\{(x,y)\in R^{* 2}|v(x)>v(y)\},
\end{eqnarray*}
then $X_1,X_2$ form a partition of $R^{* 2}$. The isomorphisms
\begin{eqnarray*}
\{0\}\times R^{* 2}\to X_1: (0,x,y)\mapsto (x,xy),\\
\{1\}\times R^{* 2}\to X_2: (1,x,y)\mapsto (\pi xy,y),
\end{eqnarray*}
imply that $R^{* 2}\sqcup R^{* 2}$ is isomorphic to $X_1\cup X_2=R^{* 2}$.
It follows that $2[R^{* 2}]=[R^{* 2}]$, so $[R^{* 2}]=0$. Notice that the
proof of (ii) does not use the full power of ($*$), only that $R$ is
definable.
\end{proof}

\begin{theorem}\label{Gring=0}
Let $K$ be a $\Z$-valued field, which is a model for the language $\Lm$.
If the structure $(K,\Lm)$ satisfies condition ($*$), then the
Grothen\-dieck ring $K_0(K)$ is trivial and there exists an isomorphism
from $R^2\setminus\{(0,0)\}$ to $R^2$.
\end{theorem}

\begin{proof} Since $0=[R]=[R^*]+[\{0\}]$ we have $[R^*]=-1$. Together with
$0=[R^{* 2}]=[R^*]^2$ this yields $1=0$, so $K_0(K)$ is trivial.

Define the isomorphisms $\psi:R^2\to \pi^3 R^2:(x,y)\mapsto (\pi^3
x,\pi^3 y)$ and $\varphi_i:R^2\to (\pi^i+\pi^3 R)\times(\pi^i+\pi^3
R):(x,y)\mapsto (\pi^i+\pi^3 x,\pi^i+\pi^3 y)$ for $i=1,2$.

Since clearly $\psi(R^*\times R^*)\cup \varphi_1(R^*\times R^*)$ is
isomorphic to $R^{*2}\sqcup R^{*2}$, we can find by
Proposition~\ref{bijections}(ii) an isomorphism
\[
f_1:\varphi_1(R^*\times R^*)\to \psi(R^*\times R^*)\cup
\varphi_1(R^*\times R^*).\] Define $f_2$ by
\[
f_2:\psi(R\times R^*)\cup\varphi_2(R^{(1)}\times R^*) \to
\varphi_2(R^{(1)}\times R^*):
\]
\[
 \left\{\begin{array}{rcl}
              \psi(x,y) & \mapsto & \varphi_2(1+\pi x,y),\\
              \varphi_2(x,y) & \mapsto & \varphi_2(\pi x,y).\end{array}\right.\
\]
Analogously, we can modify the function given in the proof of
Proposition~\ref{bijections}(i) to get an isomorphism
\[
f_3:\varphi_2(\{0\}\times R^{(1)})\to\varphi_2 (\{0\}\times
R^{(1)})\cup\psi(\{0\}\times R).
\]
Finally,
\[
  g:R^2\setminus\{(0,0)\}\to R^2:x\mapsto\left\{
\begin{array}{ll}
 f_1(x) & \mbox{if }x\in \varphi_1(R^*\times R^*),\\
 f_2(x) & \mbox{if }x\in \psi(R\times R^*)\cup\varphi_2(R^{(1)}\times R^*),\\
 f_3(x) & \mbox{if }x\in \varphi_2(\{0\}\times R^{(1)}),\\
 x      & \mbox{else},\\
\end{array}\right.
\]
is the required isomorphism.
\end{proof}
We give some examples for the conditions of Theorem~\ref{Gring=0} to be
satisfied. Let $\Lm_{ac}$ be the language of rings with an extra constant
symbol to denote $\pi$ and a relation symbol to denote the set $R^{(1)}$.
Let $\Lm_{ac,R}$ be the language $\Lm_{ac}$ with an extra relation symbol
to denote $R$.
\begin{itemize}
\item Let $K$ be a valued field with valuation to the integers
$\Z$. Then we can define an angular component as follows. Choose
$\pi\in K$ with $v(\pi)=1$ and put $ac(x)=\overline{\pi^{-v(x)}x}$
for $x\not=0$. Then clearly $ac(\pi)=1$ and $(K,\Lm_{ac,R})$
satisfies condition ($*$). \item Let $K$ be a Henselian field with
valuation to the integers $\Z$. Then $R$ is already definable in
the language of rings: if char$(\bar{K})\not=2$ we have $R=\{x\in
K|\exists y\in K,y^2=1+\pi x^2\}$ and if char$(\bar{K})=2$ then we
use the formula $\exists y\in K,y^3=1+\pi x^3$ to define $R$. This
implies that $(K,\Lm_{ac})$ satisfies condition ($*$). \item For
definability of the valuation ring in fields of rational functions
within the language of rings, see \cite{D} and \cite{KR}.
\end{itemize}

Now we specialize our attention to local fields with finite residue field.
\begin{theorem}\label{FF}
Let $K=\F_q((t))$ be the formal Laurent series over the finite field
$\F_q$ and $\Lm_t$ the language of rings with a constant symbol to denote
$t$. Then $K_0(K)$ is trivial and we have an isomorphism $R\to R^*$.
\end{theorem}
\begin{proof}
We first show that $K$ satisfies condition ($*$). Since $K$ is a Henselian
field, $R$ is definable as shown above. For each $x\in\F_q$ we have
$x^{q-1}=1$, so we can define $R^{(1)}$ as
\[R^{(1)}=\{x\in R|\exists y\in R^*,\ \bigvee_{n=0}^{q-2} t^ny^{q-1}=x\},\]
 again by Hensel's lemma.

By Theorem~\ref{bijections} we have an isomorphism $f:R^2\to
R^2\setminus\{(0,0)\}$. For a Laurent series $H(t)\in K$ we have
$H(t)^p=H(t^p)$. Consequently, the map
\[g:K^2\to K:(x,y)\mapsto x^p+ty^p\]
is an injection from the plane into the line. We obtain the isomorphism
\[R\to R^*:x\mapsto \left\{
\begin{array}{ll}
 g\circ f\circ g^{-1}(x) & \mbox{if } x\in g(R^2),\\
 x & \mbox{else.}
\end{array}\right.
\]
\end{proof}

Now let $\Q_p$ be the field of $p$-adic numbers and $K$ a fixed finite
field extension of $\Q_p$. Choose an element $\pi$ with $v(\pi)=1$, then
$ac(x)=\pi^{-v(x)}x\ \mathrm{ mod }(\pi)$  defines an angular component for
$x\not=0$. We work with $\Lm_\pi$, the language of rings with an extra
constant symbol to denote $\pi$. For a definable set $X\subset K$ and
$k\in\N_0$ we write
\begin{eqnarray*}
X^{(k)} & = & \{x\in X|x\not=0 \mbox{ and } v(\pi^{-v(x)}x-1)\ge k\},
\end{eqnarray*}
which corresponds with our previous definition of $R^{(1)}$. The set $R$
and each $X^{(k)}$ is definable by the same argument as in the proof of
Theorem \ref{FF}, so $(K,\Lm_\pi)$ satisfies condition $(*)$. We put
$P_n=\{x\in K^\times|\exists y\in K,\ y^n=x\}$ and $\bar P_n=P_n\cap R$.
Recall that $P_n$ is a subgroup of finite index in $K^\times$ for each
$n$.

For convenience, we recall the following easy corollary of Hensel's Lemma.
\begin{cor}\label{corhensel}
Let $n>1$ be a natural number. For each $k>v(n)$, and $k'=k+ v(n)$ the
function
$$
K^{(k)}  \to P_n^{(k')}:x\mapsto x^n
$$
is an isomorphism.
\end{cor}

In the next proposition we exhibit some isomorphisms between definable
sets.
\begin{proposition}\label{p-adic} Let $K$ be a finite field extension of
the $p$-adic numbers and $\Lm_\pi$ the language of rings with an extra
constant symbol to denote $\pi$. Then we have
\item{(i)} for each $k>0$, the union of two disjoint
copies of $R^{(k)}$ is isomorphic to $R^{(k)}$;
\item{(ii)} the union of two disjoint
copies of $R^*$ is isomorphic to $R^*$.
\end{proposition}
\begin{proof}
(i) \textbf{Case 1: $p\not=2$.} The map $R^{(k)}\to\bar
P_2^{(k)}:x\mapsto x^2$ is an isomorphism for each $k>0$ by Corollary
\ref{corhensel}. By Hensel's Lemma, $R^{(k)}=\bar P_2^{(k)}\cup\pi\bar
P_2^{(k)}$ is a partition. Hence the function
\[
 \{0\}\times R^{(k)}\cup \{1\}\times R^{(k)}\to R^{(k)}:\left\{\begin{array}{rcl}
                    (0,x) & \mapsto & x^2,\\
                    (1,x) & \mapsto & \pi x^2,\end{array}\right.\
\]
is an isomorphism.

\textbf{Case 2: $p=2$.} The map $R^{(k)}\to \bar P_3^{(k)}:x\mapsto x^3$
is an isomorphism by Corollary \ref{corhensel}, and by Hensel's Lemma
$R^{(k)}=\bar P_3^{(k)}\cup\pi\bar P_3^{(k)}\cup\pi^2\bar P_3^{(k)}$ is a
partition. Explicitly, we see that cubing and multiplying by $1$, $\pi$
or $\pi^2$ is an isomorphism from three disjoint copies of $R^{(k)}$ to
$R^{(k)}$.  First suppose that $k>v(2)$ and put $k'=k+v(2)$, then
$R^{(k)}\to\bar P_2^{(k')}:x\mapsto x^2$ is an isomorphism  by Corollary
\ref{corhensel}. By Hensel's lemma, we have a partition
$R^{(k)}=\bigcup_{i=1}^{2^l}\alpha_i\bar P_2^{(k')}$ for some $l\in\N_0$.
Thus we can say there are isomorphisms from $R^{(k)}$ to $2^l$ disjoint
copies of $R^{(k)}$ and to three disjoint copies of $R^{(k)}$. Some
arithmetic on the number of disjoint copies yields the required
isomorphism for $k>v(2)$.

 If $k\leq v(2)$ then $R^{(k)}$ admits a finite partition into parts of the form
$\alpha R^{(v(2)+1)}$, with $v(\alpha)=0$, and hence that the required
isomorphism exists follows from property (i) for $R^{(v(2)+1)}$.

(ii) Since  $R^*$ admits a finite partition with parts of the form
$\alpha R^{(1)}$ with $v(\alpha)=0$, this follows from (i).
\end{proof}
Now we give the solution of the problems raised by J. Denef and L.
B\'elair.
\begin{theorem}
Let $K$ be a finite field extension of $\Q_p$ and $\Lm_\pi$ the language
of rings with an extra constant symbol to denote $\pi$. Then $K_0(K)=0$
and we have an isomorphism from $R$ to itself minus a point.
\end{theorem}
\begin{proof}
The triviality of the Grothendieck ring follows from
Theorem~\ref{Gring=0}.

 We write the isomorphism explicitly in the case
$p\ne 2$. First let
\[ W = 1+ \pi^2R^* \cup \pi^2R \cup \pi + \pi^2R^{(1)}.
\]
As in the proof of Proposition~\ref{p-adic}, we can write
\[ R^* =\bigcup_{i=1}^l \alpha_iR^{(1)} =
       \bigcup_{i=1}^l\big( \alpha_i{\bar P_2^{(1)}}
                             \cup \pi\alpha_i{\bar P_2^{(1)}}\big)
\]
as a partition for some $l\in \N_0$. Thus the function
\[
   f_1:\pi^2 R^* \cup 1+\pi^2 R^*\to 1+\pi^2 R^*:
                          \left\{\begin{array}{rcl}
                 \pi^2\alpha_i x & \mapsto & 1+\pi^2(\alpha_i x^2),\\
                 1+\pi^2\alpha_i x & \mapsto & 1+\pi^2(\pi\alpha_i x^2),
\end{array}\right.
\]
where $x\in R^{(1)}$, is a well-defined isomorphism. Modify the function
given in the proof of Proposition~\ref{bijections}(i), to get
\[
  f_2:\pi^2 R\cup \pi+\pi^2 R^{(1)}\to \pi+\pi^2 R^{(1)}:
                       \left\{\begin{array}{rcl}
               \pi^2 x & \mapsto & \pi+\pi^2(1+\pi x),\\
               \pi+\pi^2 x & \mapsto & \pi+\pi^2(\pi x).
\end{array}\right.
\]
Then the function
\[
  f:W\to W\setminus\{0\}:x\mapsto\left\{
                          \begin{array}{ll}
           f_1^{-1}(x) & \mbox{if }x\in 1+\pi^2 R^*,\\
           f_2(x) & \mbox{if }x\in\pi^2 R\cup \pi+\pi^2 R^{(1)},\\
                          \end{array}\right.
\]
is an isomorphism. Finally,
\[
  g:R\to R^*:x\mapsto\left\{
\begin{array}{ll}
 f(x) & \mbox{if }x\in W\\
 x      & \mbox{if } x\notin W\\
\end{array}\right.
\]
is an isomorphism.

In the case $p=2$, we know from Proposition \ref{p-adic}(i) that there is
a function which plays the role of $f_1$. The rest is as above.
\end{proof}
\begin{remark}
\begin{itemize}
\item
The construction of the bijection $\Z_p\to\Z_p\setminus\{0\}$ also works
for the field $\F_q((t))$ if $2\nmid q$. If $2|q$ the proof of
Proposition~\ref{p-adic}(i) collapses since the index of the squares in
$\F_q((t))^\times$ is infinite.
\item The triviality of the Grothendieck ring of a structure $\M$
implies that every  Euler characteristic on the definable sets is trivial.
An Euler characteristic is a map $\chi:\Def(\M)\to R_\chi$ with  $R_\chi$
a ring, such that $\chi(X)=\chi(Y)$ if $X\cong Y$, $\chi(X\cup
Y)=\chi(X)+\chi(Y)$ if $X\cap Y=\phi$ and $\chi(X\times
Y)=\chi(X)\chi(Y)$. In general an  Euler characteristic on $\Def(\M)$
factorizes through $\Def(\M)\to K_0(\M):X\mapsto [X]$.
\end{itemize}
\end{remark}

\bibliographystyle{asl}
\bibliography{presbib}

\end{document}